\input amssym.def
\input amssym.tex
\magnification = 1200
\topskip = 10 truemm

\font \titre = cmbx10 scaled 1400

\font\bfdouze=cmbx12

\centerline{\titre  COMPACT K\"AHLER THREEFOLDS}

\

\centerline{\titre OF $\pi_1$-GENERAL TYPE.}

\bigskip

\centerline{\bfdouze F. CAMPANA and Q.ZHANG}

\vskip 1cm

 {\bfdouze 0 - Introduction.}

\medskip

 Let X be a compact complex manifold, connected and of complex dimension
n. Let $\kappa (X)$ be its Kodaira dimension, and $\alpha_X:X\to Alb(X)$ its Albanese map.
A fibration $f:X\to Y$ is a surjective meromorphic map with connected fibres.

Let the universal cover of $X$ be $\widetilde X$. We say that X is of 
{\bf $\pi_1$-general type} if 
$\widetilde X$ is not covered by its positive-dimensional compact analytic subsets. 

Obvious examples of such $X's$ are complex tori and compact quotients of bounded domains in 
$\Bbb C^n$. We conjecture that the general case is essentially a simple combination of these:

{\bf Conjecture 0:} {\it Let $X$ be a compact K\"ahler manifold of $\pi_1$-general type. Then:

1. $\kappa(X)\geq 0$.

2. There exists a finite \'etale cover of $X$ which is bimeromorphic to a manifold $X'$ such 
that its Iitaka-Moishezon fibration (usually called the Iitaka fibration. See however [Moi])
 $f:X'\to Y'$ is a submersion with fibres complex tori, and with 
base $Y'$ a projective manifold of general type and of $\pi_1$-general type.}

\

This conjecture extends the one made in [C95]: {\it if $X$ is of $\pi_1$-general type, and if 
$\kappa(X)=0$, then $X$ is bimeromorphically covered by a torus.} Due to [K-O], 
the above conjecture also extends the one by S.Iitaka ([I]: {\it if $\widetilde X\cong \Bbb C^n$, 
then $X$ is covered by a torus).}

See \S 4. below for a motivation of this conjecture, and of the terminology.

Notice that the above conjecture already appears in [K] for $X$ projective.

For curves , the conjecture holds by the
uniformisation theorem. The conjecture can easily be shown for
surfaces as well, using the Enriques-Kodaira-Shafarevitch classification.

We show here that the conjecture holds for threefolds also. 

More precisely:

{\bf Theorem 0:} {\it The above conjecture holds in any of the following three situations:

1. $X$ is projective, $K_X$ is nef, and the Abundance conjecture holds.

2. $X$ is a projective threefold.

3. $X$ is a compact non-projective K\"ahler threefold.
}

In the projective case, we just have to tie up together results already well-known, among which
 [K,6.3] (which solves the case where the generic fibre of $f$ is an Abelian variety) plays a crucial role. 
For non-projective 
threefolds, we need to use the results of [C95] to rule out the (conjecturally non-existing) 
simple non-Kummer case. We indeed show in arbitrary dimension 
that $\chi (\cal O$$_X)\neq 0$, if $X$ is not of general type, but of $\pi_1$-general type, 
which implies 
in dimension three that either $p_g:=h^{0,3}\geq 2$, or $q:=h^{1,0}\geq 1$.

The results of [CP] immediately rule out the case $\kappa(X)=-\infty$, and solve the case 
$\kappa(X)=0$. To deal with the remaining cases $\kappa(X)=1,2$, we study the 
Iitaka-Moishezon fibration. One can reduce to the case where its generic fibre 
is a complex torus, and adapt most of the arguments of [K,6.3]. We actually obtain 
a more precise result, in this case:

\

{\bf Theorem 1:} {\it Let $X$ be a compact 
K\"ahler threefold of $\pi_1$-general type. Then: 

1. $\kappa(X)\geq 0$, and a suitable finite \'etale cover $X'$ of $X$ has $f':X'\to Y'$, its Iitaka-Moishezon fibration, bimeromorphic:

2. Either: to $Alb(X')\times_{Alb(Y')}Y'$, and $Y'$ is of general type and of $\pi_1$-general type.

3. Or: to a  submersion with all fibres Abelian surfaces, and base $Y'$ a curve of genus $2$ or more.

} 

\

Notice that the proof uses properties very specific to the three-dimensional case. 

\

From this and the arguments of the proof (especially the proof of 5.5), 
we immediately deduce a concrete description of such threefolds (see \S 6).

\

{\bf Corollary 2:} {\it Let $X$ be a compact K\"ahler threefold of $\pi_1$-general type. 
Then, after a finite \'etale cover, $X$ is bimeromorphic to one of the following:

1. A complex torus.

2. Either:

2.a. A subvariety of Kodaira dimension $1$ of a complex torus, or:

2.b. A threefold with a submersion onto a curve of general type with Abelian surfaces as fibres.

3. A variety $V$ of type 3.a, 3.b, or 3.c below:

3.a. $V$ has a finite map onto a threefold of Kodaira dimension $2$ of a complex torus.

3.b. $V=E\times S$, where $E$ is an elliptic curve, and $S$ a surface of general type 
and $\pi_1$-general type such that $q(S')=0$, for any finite \'etale cover $S'$ of $S$.

3.c. $V=S\times _{C} S'$, where $S,S'$ are two surfaces fibered over a curve $C$ of general type, with $S$ of 
general type and $\pi_1$-general type, and $S'$ a surface in a complex torus, elliptic over $C$.

4. A variety of general type. }

\

In particular, S.Iitaka's conjecture is true for compact K\" ahler threefolds (this is shown also in [N"], some intermediate results of which are similar to some of ours).

We now prove Theorems 0 and 1, according to the three cases of theorems 0.

\

 {\bf \S1 - $X$ PROJECTIVE, $K_X$ NEF, ASSUMING ABUNDANCE.}

\

Because $X$ is of $\pi_1$-general type (or because $K_X$ is nef), 
it is not uniruled. In particular, if the abundance
conjecture holds for $X$, then the Kodaira dimension $\kappa X$ of $X$
should be nonnegative.

First, let us assume that $\kappa (X)=0$, because $K_X$ is nef, from Abundance we deduce that 
$K_X$ is torsion. Calabi-Yau and a result of
Beauville and Bogomolov [B] tell us that (after a finite \'etale cover), $X$ will split as
the product of some manifold with a simply-connected Calabi-Yau
manifold. Because $X$ is of $\pi_1$-general type, the second factor is a point, 
and $X$ is covered by a complex torus.

If $X$ is of general type, there is nothing to prove.

Now, let us assume that: $n>\kappa X>0$.  Then
Abundance conjecture assures that
there exists a morphism (the Iitaka-Moishezon fibration) $f:
X\rightarrow Y$. Its generic 
fibre is smooth with torsion canonical bundle. 
By the previous step, this generic fibre is an Abelian variety. 
From [K, 6.3], we get the conclusion.

\bigskip

{\bf \S 2. PROJECTIVE THREEFOLDS.}

\

The above arguments apply in this case, because the Abundance conjecture is known in this case, 
by the work of the Japanese School (Kawamata, Miyaoka, Mori) and V.Shokurov. 
We just need to deal 
with the case $\kappa(X)=0$, $X$ of $\pi_1$-general type. Then [N] shows that 
a finite \'etale cover of $X$ is bimeromorphic to a torus.

\bigskip

{\bf \S 3. K\" AHLER THREEFOLDS.}

\bigskip

Let now $X$ be a compact K\" ahler (connected smooth)  
threefold of $\pi_1$-general type with $\kappa(X)\in \lbrace -\infty, 0,1,2\rbrace$.

We now make a case-by-case study, according to the value of $\kappa(X)$. We list in this \S 3. the 
main intermediate steps of the proof, leaving the proofs of the technical 
results for later sections.

{\bf 3.1: Case $\kappa(X)=-\infty$.} 

\
Then $X$ is not projective, by \S 2.
This case is also excluded in the non-projective case by [CP,8.1]. Indeed, by loc.cit. we get that 
$X$ is either uniruled or ``simple" (which means: not covered by compact analytic subsets of 
intermediate dimensions either $1$ or $2$). In particular, if $X$ is simple,
 the algebraic dimension $a(X)$ vanishes 
(ie: all global meromorphic functions on $X$ are constant), and there is no surjective 
meromorphic fibration $f:X\to Y$ with $Y$ a curve or a surface. This ``simple" case 
is excluded by the following lemma, since we assumed $\kappa(X)=-\infty$:

\

{\bf 3.1.1 Lemma:} {\it Let $X$ be a simple compact K\"ahler threefold. If $X$ is of
$\pi_1$-general type, then $X$ is bimeromorphic to a torus. }

\

{\bf Proof:} From theorem 3.1.2 below, we learn that $\chi (\cal O$$_X)=0$. Because 
$h^{2,0}\geq 1$, by the result of Kodaira recalled above, we see that $q+p_g\geq 2$, 
and so either $p_g:=h^{0,3}(X)\geq 2$, or $q:=h^{1,0}\geq 1$. The first case is excluded because 
it would produce a nonconstant meromorphic function on $X$. In the second case, the Albanese 
map $\alpha_X:X\to Alb(X)$ needs to be surjective and connected (because $a(X)=0$; see [U]). 
Because $X$ is moreover simple, its fibres are zero-dimensional (because $q=dim(Alb(X))>0$). 
Thus $\alpha_X$ is bimeromorphic. $\square$

\

We have used the following result, whose proof is given in 4.3.4):

\

{\bf 3.1.2 Theorem:} {\it Let $X$ be a compact connected K\"ahler $n$-fold of $\pi_1$-general type. 
Then: either $X$ is projective and of weak general type (ie: $\kappa^+(X)=n$, see \S 4.), or: $\chi (\cal O$$_X)=0$.}

\

{\bf 3.2 Case $\kappa(X)=0$.}

\

In the projective case, the conclusion follows from \S 2. 
This case is easily solved again using [CP, 8.1] in the non-projective case: $X$ has then to be either 
simple, or after a finite \'etale cover, bimeromorphic to either a torus, or the product of 
an elliptic curve with a $K3$ surface. Since $X$ is of $\pi_1$-general type, $X$ is bimeromorphic 
to a torus.

\

{\bf 3.3 Reduction steps when $\kappa(X)\geq 1$.} 

\

We then assume that $\kappa(X)=1$ or $2$. Let $f:X\to Y$ be the (holomorphic, 
after resolution of 
its indeterminacies) Iitaka-Moishezon fibration of $X$. Its base $Y$ is 
either a curve or a surface. Its generic 
fibre $X_y$ has $\kappa(X_y)=0$, and is of $\pi_1$-general type. Thus $X_y$ is an elliptic curve  
if $\kappa(X)=2$, and bimeromorphic to either a torus or a bielliptic surface if $\kappa(X)=1$. 
In this last case, by a finite \'etale cover of $X$, we can assume $X_y$ to be bimeromorphic 
to a torus, by using a very special case of [K, 6.4.1] . We thank N.Nakayama who pointed out a gap in the proof given in the previous version of the present text, where a more general version was stated. The special case below is the same one as in [N", 6.4]. 

\

{\bf 3.3.0. Proposition:} {\it  Let $H<G$ be groups, with $H$ normal in $G$, such that $H$ is finitely generated and residually finite. Let  $H_1$ be an abelian characteristic subgroup of finite index in $H$. Then: 

1. $G$ has a normal subgroup $G'$ of finite index such that $H':=G'\cap H\subset Z.H_1$, where 
$Z\subset H$ is the centraliser of $H_1$ in $H$ (ie: the subgroup of $H$ consisting of all elements commuting with each element of $H_1$). Moreover, $G'$ can also be chosen such that the quotient group $(G'\cap Z.H_1)/(G'\cap H_1)$ is abelian.

2. If $H_1$ is abelian and torsionfree, then $G'$ can be choosen such that $(G'\cap H)$ is abelian.}

\

{\bf Proof:}  Of (1): Let $F:=H/H_1$, and $c:H\to Out(H_1):=Aut(H_1)/Int(H_1)$ be the action by conjugation. Its kernel is $Z.H_1$, and its image $F'$ is a quotient of $F$, hence finite, because $c$ is trivial on $H_1$.

For $K$ a characteristic subgroup of finite index of $H_1$, define, by using conjugation again:
$c_K:H\to Aut(H_1/K)$. Because $H$ is residually finite, and $H_1$ is abelian, we can choose $K$ such that the kernel of $c_K$ is included in the kernel of $c$. 

Indeed, if $h\in H$ does not belong to the kernel of $c$, then $[h,a]:=hah^{-1}a^{-1}\neq 1$, for some $a\in H_1$. Thus $[h,a]\notin K$, for some suitable $K$, characteristic and of finite index in $H_1$, by the residual finiteness of $H$. And so, $h$ does not act trivially on $H_1/K'$, for any smaller $K'\subset K$, as desired. Because $F'$ is finite, we need to check this condition for finitely many $h's$ in $H$, say for $h_i,i\in H/Z.H_1$, representatives of $H/Z.H_1$, since if $h=h_i.z.b$, for $z\in Z,b\in H_1$, we have:
$[h, a]=[h_i,a].[b,a]=[h_i,a]\neq 1$, because $z\in Z$, and $H_1$ is abelian. Hence the claim.

We choose and fix $K$ such that the preceding condition is satisfied.

Define now $r:G\to Aut(H_1/K))$ by conjugation. Because the image group is finite, the kernel $G"$ of $r$ is of finite index in $G$. And $G"\cap H\subset Z.H_1$, by the above construction and observation. To ensure the second condition, we simply need to ask $G'$ to be the kernel of the conjugation map: $r':G"\to Aut(H/H_1)$.

\

 Proof of (2):  this follows from the following standard claim (applied to $A=(G'\cap H_1)$, and $B:=G'\cap Z.H_1$):
 
 {\bf Claim:} Let $A\subset B$ be groups, with $A$ central (hence abelian),  torsionfree of finite index in $B$, and such that $B/A$ is abelian. Then $B$ is abelian.
 
 {\bf Proof:} For $x,y,z \in B$, we have :$xyx^{-1}y^{-1}:=[x,y]=[xA,yA]$, and by easy commutator calculus: $[x,yz]=[x,y][x,z]^y=[x,y][x,z]$ (since $B$ acts trivially on $A$). Thus, for any $x\in B$,  $[x,B]$ is a finite subgroup of $A$. And is thus trivial, by the torsionfreeness of  $A$. $\square$

\

This being done, we can use [C85] (the construction of 
relative Albanese maps) to further reduce 
to the special case when $X_y$ is a torus.

We can and shall thus assume in the sequel that $X_y$ is a torus for $y\in Y$ generic.

\

{\bf 3.3.1 Lemma:} {\it Assume that $\kappa(X)=1$, for $X$ a compact K\" ahler threefold 
of $\pi_1$-general type. Let $f:X\to Y$ be its Iitaka-Moishezon fibration. Then $Y$ is a curve 
of general type (ie: $g(Y)\geq 2$), after $X$ has been 
replaced by a suitable finite \'etale cover.

Hence we can (and shall) assume that $f$ has no 
multiple fibres, and that $Y$ is a curve of general type, if $\kappa(X)=1$.}

\

{\bf Proof:} Assume the contrary. After replacing $X$ by a suitable 
finite \'etale cover, we can assume that $f$ is $\pi_1$-exact, which means that the 
natural sequence of groups (3.3.1.1) below is exact: 

(3.3.1.1) $\pi_1(X_y)\to \pi_1(X)\to \pi_1(Y)\to 1$.

This property is shown in [C98], and is peculiar to the case of $1$-dimensional bases.
Thus $\pi_1(X)$ is solvable of finite rank, and the Albanese map of $X$ and any of its 
finite \'etale covers is obviously surjective (because no such cover has a surjective 
map onto a positive dimensional variety of general type). In this case, it follows from [C01] and 
[C95'] that $\pi_1(X)$ is almost abelian (ie: has a subgroup of finite index which is abelian).

By taking a suitable finite \'etale cover of $X$, we shall 
thus assume $\pi_1(X)$ to be abelian, and moreover that
 $\alpha_X:X\to Alb(X)$, the Albanese map of $X$ induces an 
isomorphism of fundamental groups. Because $X$ is of $\pi_1$-general type, 
$\alpha_X$ has to be generically finite. Because the fibres of $f$ are tori, the restriction 
of $\alpha_X$ to $X_y$ has to be \'etale onto its image, and so isomorphic, because $\alpha_X$ 
induces an isomorphism of fundamental groups. Thus $\alpha_X$ is bimeromorphic, and $\kappa(X)=0$, 
in contradiction with our hypothesis. $\square$

\

{\bf 3.3.2 Proposition:} {\it Let $X$ be a compact k\"ahler 
threefold of $\pi_1$-general type, with 
$\kappa(X)\geq 1$. Let $f:X\to Y$ be its Iitaka-Moishezon fibration.

Replacing $X$ by a suitable finite \'etale cover, 
$X$ is bimeromorphic to either:

 $Alb(X)\times _{Alb(Y)} Y$, or such that $f$ is a submersion onto the 
curve of general type $Y$ with Abelian surfaces as fibres.}

\

 Proposition 3.3.2, which is the main assertion of theorem 1,
 will be proved in \S 5.

The proof of theorem 1 is then completed by means of the easy:

\

{\bf 3.3.3 Lemma:} {\it Let $f:X\to Y$ be the Iitaka-Moishezon fibration of 
$X$ k\"ahler compact. Assume that 
$X$ is bimeromorphic to $Alb(X)\times _{Alb(Y)} Y$. Then:

1. $\kappa(X)=\kappa(Y)$.

2. $X$ is of $\pi_1$-general type if and only if so is $Y$. More generally:

3. $dim(X)-\gamma d(X)=dim(Y)-\gamma d(Y)$ (see \S 4. for this notion).}

\

{\bf Proof:} Assertion 1 follows from the equality:

$K_{X/Y}=\alpha_X^*(K_{[Alb(X)/Alb(Y)]})=\alpha_X^*(\cal O$$_{Alb(X)})$.

For assertion 3, we have a natural map $h$:

$h:\Gamma(X)\to 
Z:=\Gamma(Alb(X))\times _{\Gamma(Alb(Y))}\Gamma(Y)=Alb(X)\times _{Alb(Y)}\Gamma(Y)$, 

(induced from the equality: 
$X=Alb(X)\times _{Alb(Y)}Y$). 

We shall show that $h$ is generically finite. Let $H\subset X$ be an irreducible compact 
analytic subset contained in a general fibre of the composite map $h\circ \gamma_X:X\to Z$. 
We have to show that its image in $\Gamma(X)$ is point. But this is immediate, because $H$ is 
of the form:$a\times G$, for $G\subset Y$ a general fibre of $\gamma_Y$. 
The image of $\pi_1(\hat H)$ 
in $\pi_1(X)$ is thus finite. This shows assertion 3, of which assertion 2 is a special case. $\square$

\

\

 {\bf \S 4 - VANISHING OF  $\chi (\cal O$$_{X})$.}

\

We recall here the two bimeromorphic invariants $\kappa^+(X)$ and $\gamma d (X)$, introduced 
in [C95] for compact $n$-dimensional K\"ahler manifolds $X$.

\

{\bf 4.1 The invariant $\gamma d(X)$.}

\

Let $r:\widetilde X\to X$ be the universal cover. There exists a unique proper connected 
fibration $\tilde \gamma:\widetilde X\to \widetilde \Gamma(X)$ such that if $\tilde x \in 
\widetilde X$ is a general point, any connected compact analytic 
subset of $\widetilde X$ going through $\tilde x$ is contained in the fibre of $\tilde 
\gamma$ through $\tilde x$. Because $\tilde \gamma$ is obviously $\pi_1(X)$-equivariant, 
the map $\tilde \gamma$ induces a fibration,called the $\Gamma$-reduction of $X$:
$\gamma_X:X:=(\widetilde X/\pi_1(X))\to \Gamma(X):=\tilde \Gamma(X)/\pi_1(X)$,

 (This last map is constructed 
independently, for $X$ projective, in [K] and termed: ``Shafarevitch map of $X$").

The geometric interpretation of $\gamma$ is that if $V\subset X$ is a compact 
irreducible analytic subset going through a general point $x\in X$, and 
$\hat V$ its normalisation, then $\pi_1(V)_X\subset \pi_1(X)$, the image of 
the natural map: $\pi_1(\hat V)\to \pi_1(X)$, is finite if and only 
if $V$ is contained in the fibre of $\gamma$ through $x$.

The fibration $f:X\to Y$ naturally induces a natural map $\gamma_f:\Gamma (X)\to \Gamma(Y)$.

\

{\bf 4.1.1 Definition:} We denote by $\gamma d(X):=dim(\Gamma(X))\in\lbrace 0,1,...,n\rbrace$, 
and call it the $\Gamma$-dimension of $X$.

\

Thus $\gamma d(X)=0$ iff $\pi_1(X)$ is finite, 
and $\gamma d(X)=n$ iff $X$ is of $\pi_1$-general type.

Thus $(dim(X)-\gamma d (X))$ is the maximal dimension of subvarieties of $X$ having 
``finite fundamental group in $X$".

\

{\bf 4.2 The invariant $\kappa^+$.}

\

{\bf Definition:} For $X$ a compact complex connected manifold, we define :
$\kappa^+(X):= max\lbrace \kappa (X,Det(F)), F\subset \Omega_X^p,p>0\rbrace$, here 
$F$ runs through all coherent non-zero subsheaves of $\Omega_X^.$, $Det(F)$ is the 
saturation of its determinant in $\Lambda ^r(\Omega_X^.)$, and 
the Kodaira dimension is suitably defined (see [C95] for details).

\

We of course have: $\kappa^+(X)\geq \kappa(X)$, and conjecturally (see [C95]): 
$\kappa^+(X)=\kappa(X)$ if $\kappa(X)\geq 0$, while $\kappa^+(X)=-\infty $ iff 
$X$ is rationally connected. These conjectures are consequences of 
the standard conjectures of MMP when $X$ is projective, at least.

\

One extreme case is proved in [K], with a different formulation:

\

{\bf 4.2.3 Proposition:} {\it  Assume that $X$ is 
compact K\"ahler, $n$-dimensional, 
and such that: $\kappa^+(X)=n$, with $K_X$ nef. Then: $\kappa(X)=n$.}

\

{\bf Proof:} Let $\cal F$ $\subset $${\Omega }^{p} _{X} ,(p>0)$
be a coherent subsheaf of rank r ,
such that:

 $\kappa (X,Det(\cal F$$))= n$,where Det($\cal F$) is the saturation
of  det($\cal F$) in :
$\Lambda ^r ( {\Omega }^{p} _{X} $).

By [Mi], if $C$ is the generic curve cut out by sufficiently ample linear
systems, then
$\Omega ^1_{X}$, and so also $\Omega ^p_{X}$ restricted to $C$ are
semi-positive. The quotient (on $C$
locally free)
sheaf Q:=  $\Omega ^p_{X} / \cal F $ has thus nonegative degree on C (ie:
deg($Q _{/C}) \geq 0$).

Passing to the limit, in the very ample linear systems used, we obtain the inequality:
\medskip
(4.2.4) D.$N_1.N_2....N_{n-1} \geq 0$ for arbitrary nef divisors
$N_1,...,N_{n-1}$ ,where D is the first Chern class of Q.

\medskip
Now, Det($\cal F)$ is big by assumption, and so there exist $\Bbb Q$-divisors A and E,with
$A$ ample and $E$ effective, such
that Det($\cal F)$= $A+E$.

For a suitable positive integer m, we have:

$m.K_X=det (\Omega ^{p}_{X})= Det(\cal F$$)+D= A+E+D$.

And so: $(m.K_X )^n$-$A ^n$= (E+D).($\Sigma ^{i= n-1}_{i= 0} A ^{i}.(m.K_X)
^{n-1-i})\geq 0$, by (4.2.4).

Thus $K_X ^n \geq m ^{-n}.A ^n >0$ , as claimed. $\square$

\

{\bf 4.3 The comparison theorem:}

\

It gives an inequality between $\kappa^+$ and $\gamma d$.

\

{\bf 4.3.1 Theorem ([C95]):}{\it Let $X$ be a compact K\"ahler manifold with 
$\chi(\cal O$$_X)\neq 0$.

Then $\kappa^+(X)\geq \gamma d(X)$.}

\

Observe that if $X$ is a complex torus, we have: $\gamma d(X)=n$, and $\kappa(X)=\kappa^+(X)=0$.
Of course, $\chi(\cal O$$_X)= 0$, then. It was conjectured in [C95] that $X$ is (up to 
bimeromorphy) covered by a torus if $\gamma d(X)=n$, and $\kappa(X)=0$.

Thus $\gamma d(X)$ behaves in some sense like a Kodaira dimension as well. This motivates the term 
``of $\pi_1$-general type", and the conjecture 0 above.

\

{\bf 4.3.2 Corollary:} {\it Let $X$ be a compact K\"ahler $n$-fold. Assume that $X$ is 
of $\pi_1$-general type (ie: $\gamma d(X)=n$), and that $K_X$ is nef. Then: either 
$\kappa(X)=n$, or $\chi(\cal O$$_X)= 0$.}

\

{\bf Proof:} Assume $\chi(\cal O$$_X)\neq 0$. Then $\kappa^+(X)=n$, 
by the comparison theorem 4.4.1. In particular, $X$ is projective. Because $K_X$ is nef, 
we obtain from 4.2.3 that $\kappa(X)=n$. $\square$

\

(Observe that we used the hypothesis $K_X$ nef only in the projective case. 
So it is not needed to give the delicate definition in the K\"ahler case).
Of course, we conjecture that the above result holds without the hypothesis that $K_X$ is nef.

\

{\bf 4.3.3 Corollary:}{\it Let X be a compact K\"ahler non projective $n$-fold of 
$\pi_1$-general type. Then $\chi ({\cal O}_{X}) = 0$. }

{\bf Proof:} Assume on the contrary that $\chi (\cal O$$_{X})\neq 0 $. We have:
$\kappa ^+ (X)\geq \gamma d (X)= n$. And so $X$ is projective. Contradiction.
Thus $\chi(\cal O$$_X)= 0$.$\square$

\

{\bf 4.3.4 Example:} We shall consider here the very simple case of surfaces. Let $X$ be a 
minimal compact K\"ahler surface of $\pi_1$-general type, and of non general type 
(ie: $\kappa(X)<2$). 
We know then 
that $X$ is not uniruled. If $\kappa(X)=0$, then $X$ is covered by a torus, by classification. 
If $\kappa(X)=1$, using $\chi(\cal O$$_X)= 0$ 
(because for surfaces, $\kappa^+=\kappa if \kappa\geq 0$) and [B-P-V, III,18.2.3 and V.12.3], 
we see that if $f:X\to Y $ is 
the Iitaka-Moishezon fibration, 
then $f_*(\omega_{X/Y})$ is torsion, and has only fibres either smooth or multiple elliptic 
(of type $_mI_0$). And so, after a finite \'etale cover of $X$, 
$f$ is a principal fibre bundle over $Y$ with 
fibre some fixed elliptic curve.

\

{\bf 5. PROOF OF 3.3.2 .}

\

{\bf 5.0 Situation.} Recall that $X$ is a compact K\"ahler threefold of 
$\pi_1$-general type, with $\kappa(X)=1$ or $2$, 
and that its Iitaka-Moishezon fibration $f:X\to Y$ 
either has base a curve $Y$ of general type, has no multiple fibre, and has a generic fibre which is 
a two-dimensional complex torus, or is an elliptic fibration over a surface $Y$.

We want to show that $X$ is bimeromorphic to either $Alb(X)\times_{Alb(Y)}Y$, or to a submersion 
with Abelian surfaces as fibres, after a suitable 
finite \'etale cover of $X$.

\

{\bf 5.1 Monodromies.} We shall always denote by $Y^*$ the Zariski 
open subset of $Y$ over which $f$ is smooth, and 
by $\rho:\pi_1(Y^*)\to Aut(H_1(X_{y_0}, \Bbb Z))$ the associated monodromy representation, 
if $Y_0\in Y^*$ is given. Let $G$ be its image.

\

{\bf 5.2 Local monodromies.} If $D_i$ is a component of the divisor $(Y-Y^*)$, 
and if $\gamma _i$ is a small loop winding once around $D_i$ in $Y$, we shall denote by $T_i$ (as in [K,6.3]), 
the (well-defined up to conjugation) image of $\gamma_i$ by $\rho$. We also denote by $r_i$ the multiplicity 
of the fibre of $f$ above a generic point of $D_i$.

Let $H:=\pi_1(X_{y_0})_X/Torsion$ be the image in $\pi_1(X)$ of $\pi_1(X_{y_0})$, modulo its torsion 
subgroup. Let $q:H_1(X_{y_0}, \Bbb Z)\to H$ be the natural quotient map, and $W$ its kernel.

Easy arguments (see [K,6.3] for details) show that $N:=\sum _{i\in I}Ker (1-T_i^{r_i})\subset W$

We write: $N_{\Bbb Q}:=N\otimes _{\Bbb Z} \Bbb Q$.

\

{\bf 5.3 Proposition:} {\it For $f$ as in 5.0 above, $N_{\Bbb Q}\subset H_1(X_{y_0}, \Bbb Q)$ 
is a sub Hodge structure.

(This means here that the orthogonal of $N$ in $H^1(X_{y_0}, \Bbb R)$ 
is a complex subspace, when $H^1(X_{y_0}, \Bbb R)$ is identified with $H^0(X_{y_0}, \Omega ^1_{X_{y_0}})$ in 
the Hodge decomposition of $H^1(X_{y_0}, \Bbb C)$).}

\

This holds in a much broader context, in fact (all K\"ahler fibrations, but the proof 
seems to be written only when the fibres are projective varieties, so we must appeal to other arguments 
when the fibres are not algebraic).

\

{\bf Proof:} It rests on results stated below, but whose proofs do not depend on the statement 5.3.
We distinguish three cases:

1. $\kappa(X)=1$, and $X_y$ is not projective. Then $f$ has ``constant moduli" (ie: 5.5 holds). 
The monodromy is then (dually) a finite subgroup of $Aut_{\Bbb C} (H^0(X_{y_0}, \Omega ^1))$, 
this vector space 
being identified with $(H^1(X_{y_0}, \Bbb R))$ by means of the Hodge decomposition on $X_{y_0}$. 
See the proof of 5.6 for details.

2. $\kappa(X)=2$. The fibres are projective, here. Hence the results of Deligne [D] apply and 
give the claim. We give another, more elementary argument. 
Let $C\subset Y$ be a generic very ample curve. Let $X_C=f^*(C)$. Then $X_C$ is 
an elliptic surface of $\pi_1$-general type, with $\kappa=1$. From Example 4.3.4, we learn that after 
a finite \'etale cover, this surface is a principal fibre bundle with fibre an elliptic curve. 
The monodromy is thus finite, and complex linear in this case too (because $\pi_1(C\cap Y^*)$ 
maps surjectively to $\pi_1(Y^*)$, by Lefschetz theorem, and monodromy is compatible with restrictions).

3. $\kappa(X)=1$, and $X_y$ is projective. In this case, the fibres of $f$ are projective, and so 
the results of [D] apply. Alternatively, we can use the untwisted model of $f$ (see 5.5.1 below), 
which is projective, 
and has the same monodromy as $f$ to apply the arguments of [K,6.3]. $\square$.

\

{\bf 5.4 Corollary:} {\it Let $f:X\to Y$ be as in 5.0. Then $T_i=1, \forall i\in I$.} 

\

{\bf Proof:} If $N\neq 0$, then there is a non-zero sub Hodge structure of $H^1(X_{y_0}, \Bbb Q)$ mapped 
to $1$ in $\pi_1(X)$. But this sub Hodge structure defines a non-zero complex subtorus $A$ of $X_{y_0}$, 
whose $\pi_1$ is mapped to $1$ in $\pi_1(X)$. This is a contradiction to the 
fact that $X$ is of $\pi_1$-general type.

Thus: $(1-T_i^{r_i}), \forall i\in I$. If $\kappa(X)=1$, we have $r_i=1, \forall i$, because $f$ has no 
multiple fibres. And so 5.4 holds in this case.

Now, if $\kappa(X)=2$, we saw in the proof of 5.3 that all multiple fibres of $f$ are of type $_mI_0$, which 
directly implies that $T_i=1$. $\square$

\

{\bf 5.5 Lemma:} {\it Let $f:X\to Y$ be as in 5.0. 
There is a fixed complex torus $F$ such that $X_y \cong F$, 
for any smooth fibre $X_y$ of $f$, unless maybe when $\kappa(X)=1$ and the generic fibres $X_y$ of $f$ are 
Abelian surfaces. In this case, $f$ is bimeromorphic to a submersion with Abelian surfaces as fibres.}

\

{\bf Proof:} When $\kappa(X)=2$, this is an immediate consequence of example 4.3.4, and the proof 
of 5.3 in this case, by taking the surface $X_C=f^{-1}(C)$, for $C$ a generic very ample curve of $Y$.

When $\kappa(X)=1$, and the general fibre $X_y$ of $f$ is a non-projective torus, this follows from
 [F], or [CP, 7.1]. Note that this is stated only for $f$ the algebraic reduction 
of $X$ there, but that the proof actually works for any fibration. 

The remaining case is when $\kappa(X)=1$, and $X_y$ is an abelian surface.

\

For this, we need to consider here: 

\

{\bf 5.5.1 The untwisted Model of $f$.} For any fibration $f:X\to Y$, with $X$ compact K\"ahler, and with 
generic fibres $X_y$  complex tori, we 
introduced in [CP,6.3] the so-called ``untwisted" model $f':X'\to Y'$ of $f$. 
This fibration is caracterised up to bimeromorphic equivalence by the following two properties:

1. For $y\in Y$ generic, we have: $X_y\cong X'_y$.

2. $f'$ has a meromorphic section.

The construction is sketched in [CP, 6.3]. It was used already in [F] and [C85].

Recall from [CP,6.5] that $h^p(X,\cal O$$_X)=h^p(X',\cal O$$_{X'}), \forall p\geq 0$. 

In particular, $q(X)=q(X')$.

\

{\bf 5.5.2 Comparison of monodromies.}  
Let now, as always, denote by $Y^*$ the Zariski open subset of $Y$ over which $f$ and $f'$
are smooth, and 
by $\rho:\pi_1(Y^*)\to Aut(H_1(X_{y_0}, \Bbb Z)$ and 
$\rho':\pi_1(Y^*)\to Aut(H_1(X'_{y_0}, \Bbb Z)$ the associated monodromy representations, 
if $Y_0\in Y^*$ is given. We trivially have (identifying $X_{y_0}$ and $X'_{y_0}$):

\
 
{\bf 5.5.3 Lemma:} {\it $\rho=\rho'$.}

\

{\bf Proof:} The fibration $f'$ is obtained from $f$ over $Y^*$ simply by twisting locally 
by means of translations in the fibres (as seen from the Chow Scheme construction of $f'$). 
These translations act trivially on the cohomology. $\square$

\

We can now conclude the proof of 5.5 in the case considered ($\kappa(X)=1$, 
and $X_y$ is an abelian surface):

First: because $X_y$ is projective for $y\in Y$ generic, it follows from [C81] 
that $X'$ is projective (although $X$ is not), because $f'$ has a section. Next, we have seen that the 
local monodromies of $f$ are trivial. Thus the same holds for the local monodromies of $f'$. 
Because $f'$ has a section, we learn from [K,6.8] that after a base change over a finite \'etale cover of 
$Y$, $X'$ is bimeromorphic to an abelian group scheme over $Y$, corresponding to a 
map $\mu:Y\to \cal A$$_3$, the 
moduli of polarised Abelian surfaces with a level-three structure. 

Because the fibration $f$ has no multiple fibres, and its generic fibres are complex tori, 
it has local sections (in the analytic topology) 
over its one-dimensional base. Thus $f$ and $f'$ are locally (over $Y$) bimeromorphic. 
Hence the claim.
$\square$

\

{\bf 5.6 Lemma:} {\it Let $f:X\to Y$ be as in 5.0. The monodromy becomes trivial, 
after a finite \'etale cover of $X$, and $q(X)=q(Y)+q(F)$, unless maybe when $\kappa(X)=1$ and the generic 
fibre $X_y$ of $f$ is an Abelian surface.}

\

{\bf Proof:} The second assertion follows from the first and 5.5, using [B]. For the first assertion, we show 
that the monodromy representation $\rho$  has finite image contained (dually) 
in the complex linear automorphisms of 
$H^1(X_{y_0},\Bbb R)$. This is an easy consequence of 
the compactness of the Chow Scheme of $X$ in the K\"ahler case: 
just consider indeed the connected component $\cal C$ of the Chow Scheme of $(X\times _{Y}X)/Y$ which 
for $y$ generic in $Y$ consists of isomorphisms of $F$ with $X_y$ ($F$ 
a fixed torus isomorphic to $X_y$, the existence of 
which is the content of 5.5), 
this component containing a given isomorphism of $F$ with $X_{y_0}$. The connected 
components of the fibre of $\cal C$ over $y_0$ then consists of exactly the isomorphisms of $F$ 
with $X_{y_0}$ obtained by composing the given one with the image of $\rho$, and 
with translations of $X_{y_0}$. This shows the assertion. $\square$

\

We finally obtain 3.3.1 as follows: we saw that either $\kappa(X)=1$, and 
$f$ is bimeromorphically a submersion 
onto a curve of general type with Abelian surfaces as fibres, or we reduced to the case 
where: the global monodromy of $f$ 
is trivial, the generic fibre is a fixed torus $F$, and $q(X)=q(Y)+q(F)$. This last case we now treat.

We consider the two cases $\kappa(X)=1,2$ separately again:

If $\kappa(X)=1$, then $f$ has no multiple fibres and trivial monodromy, it is then easy to check that the 
natural map: $v:X\to V:=Alb(X)\times _{Alb(Y)} Y$ is bimeromorphic 
(the details for the case $\kappa(X)=2$ below apply to show this).

If $\kappa(X)=2$, we have to eliminate the multiple fibres of $f$ by a finite \'etale cover of $X$. In 
this case, the natural map $v:X\to Alb(X)\times _{Alb(Y)} Y$ is surjective and generically finite (by 5.6).

Moreover, its restriction to a generic fibre $X_y$ of $f$ is \'etale, and maps $X_y$ to the corresponding 
fibre $V_y$ of $V$ over $Y$. The natural map: $Alb(v):Alb(X)\to Alb(V)$ is isomorphic, by construction, 
but $\pi_1(Alb(X_y))$ maps to a subgroup $K'$ of finite index in
 $K:=Ker(\pi_1(Alb(X))\to \pi_1(Alb(Y))$. Choose a 
subgroup of finite index $A'$ in $\pi_1(Alb(X))$ such that $A'\cap K=K'$, and 
let the corresponding \'etale cover of $Alb(X)$ be $Alb'(X)$. Replacing $X$ by 
$X':=Alb'(X)\times _{Alb(X)} X$ then gives the desired \'etale cover of $X$. This is because the natural map 
from $X'$ to $V':=Alb'(X)\times _{Alb(Y')} Y$, (where $Y'$ 
is the Stein factorisation of the composite map $X'\to Y$) maps bijectively $X'_{y'}$ to $V'_{y'}$.

\

{\bf 6. PROOF OF COROLLARY 2.}

\

We treat the situation case-by-case, assuming $\kappa(X)<3$. We know from theorem 1 that $\kappa(X)\geq 0$, 
and that $X=Alb(X)\times _{Albb(Y)} Y$, if $f:X\to Y$ is the Iitaka-Moishezon fibration.

If $\kappa(X)=0$, then $Y$ is a point, and $X=Alb(X)$.

If $\kappa(X)=1$, then $Y$ is a curve of general type (by 3.3.1) isomorphic to $Alb(Y)$. 
Thus either $X=\alpha_X(X)\subset Alb(X)$, the Albanese image of $X$, or $X$ has a submersion onto the curve 
$Y$ of general type with fibres Abelian surfaces (cases 2.a, 2.b respectively).

If $\kappa(X)=2$, we shall distinguish three cases, according to $d:dim(\alpha_{Y} (Y))$, the dimension of the Albanese 
image of $Y$.

3.a occurs of course exactly when $d=2$, so that $\alpha_Y$ is generically finite.

3.b occurs when $d=0$, for any finite \'etale cover of $Y$.

3.c occurs when $d=1$, so that $C$ is the Albanese image of $Y$. (Notice that 
 $C$ is indeed of general case, after a suitable finite \'etale cover of $Y$, because $\kappa(Y)=2$). $\square$.

\

\centerline {\bfdouze Bibliographie}

\vskip 1cm

\item{\bf [B-P-V]}W.BARTH-C-PETERS-A.VAN DE VEN.{\it Compact Complex
Surfaces.Springer-Verlag 1984.}

\item{\bf [B]}A.BLANCHARD. {\it Sur les vari\'et\'es analytiques complexes. Ann.Sc.ENS 73 (1956), 157-202.}

\item{\bf [81]}F.CAMPANA. {\it Cor\'eduction alg\'ebrique d'un 
espace aalytique faiblement K\"ahl\'erien compact. Inv. Math. (1981), 187-223.}

\item{\bf [C85]} F. CAMPANA {\it R\'eduction d'Albanese d'un morphisme K\"ahl\'erien propre 
et applications, I et II.Comp. Math.54 (1985), 373-398 et 399-416.}

\item{\bf [C95]} F. CAMPANA {\it Fundamental group and positivity properties
of cotangent bundles of compact K\"ahler
manifolds. J. Alg. Geom., {\bf 4}, p. 487-502, (1995).}

\item{\bf [C95']} F. CAMPANA {\it Remarques sur les groupes de K\"ahler
nilpotents. Ann. ENS, {\bf 28}, p. 307-316, (1994).}

\item{\bf [C98']} F. CAMPANA {\it Negativity properties of compact curves in
infinite covers of projective surfaces. J. Alg. Geom.{\bf 7}, p. 673-693, (1998).}

\item{\bf [C01]} F. CAMPANA {\it Green-Lazarsfeld sets and solvable
quotients of K\"ahler groups. (10) J. Alg. Geom. (2001), 599-622.}

\item{\bf [CP]} F.CAMPANA-T.PETERNELL. {\it Complex Threefolds with non-trivial 
Holomorphic 2-Forms. J.Alg.Geom. {\bf 9}(2000), 223-264.}

\item{\bf [D]}P.DELIGNE.{\it Theorie de Hodge II. Publ.Math.IHES {\bf
40},5-57 (1972)}

\item{\bf [I]}S.IITAKA {\it Genera and classification of Algebraic Varieties. 
Sugaku 24 (1972), 14-27. }

\item{\bf [Kaw]} Y. KAWAMATA {\it Characterization of Abelian Varieties.
Comp. Math., {\bf 43}, p. 253-276, (1981).}

\item{\bf [K]} J.KOLL\`AR {\it Shafarevitch maps and plurigenera of Algebraic
Varieties.Inv.Math., {\bf 113}, p. 177-215,(1993).}

\item{\bf [K-O]}KOBAYASHI-OCHIAI.{\it Meromorphic Mappings 
into compact complex spaces of general type.Inv.Math. 31 (1975),7-16.}

\item{\bf [Mi]}Y.MIYAOKA{\it On the Kodaira Dimension of a Minimal 
Threefold. Math. Ann. 281 (1988), 325-332. }

\item {\bf [Moi]}B.MOISHEZON {\it Algebraic Varieties and Compact Complex Spaces. Actes du Congr\`es 
International des Math\'ematiciens. Nice 1970. Vol. 2, 643-648.  }

\item{\bf [Mo]}S.MORI {\it Flip Theorem and the existence 
of minimal models for threefolds. J.AMS. 1 (1988),117-253.}

\item{\bf [N']}N.NAKAYAMA.{\it Projective algebraic varieties whose universal covering spaces are biholomorphic to $\Bbb C^n$.J.Math.Soc.Japan 51 (1999),645-654.}

\item{\bf [N"]}N.NAKAYAMA.{\it Compact K\" ahler manifolds  whose universal covering spaces are biholomorphic to $\Bbb C^n$.RIMS preprint-1230(1999),1-100.}

\item{\bf [N]}Y.NAMIKAWA.{\it Deformation of Calabi-Yau Threefolds. Preprint (1993)}

\item{\bf [U]} K. UENO {\it Classification theory of Algebraic varieties
LNM, {\bf 439}, Springer-Verlag (1975).}

\vskip 1cm

{\bf e-mail}  F. CAMPANA: campana@iecn.u-nancy.fr
\vskip 1cm

{\bf e-mail}  Q.ZHANG: qzhang@math.missouri.edu
\end